\documentclass[twoside]{article}

\usepackage{amssymb}
\usepackage{amstext}
\usepackage{amscd}
\usepackage{makeidx}
\usepackage{amsmath,amsthm}
\usepackage[dvips]{graphicx}

\newtheorem{theorem}{Theorem}

\newtheorem{lemma}{Lemma}
\newtheorem{remark}{Remark}

\newtheorem{definition}{Definition}
\begin{document}

\title{Disaggregation     of      Long     Memory     Processes     on
$\mathcal{C}^{\infty}$ Class}

\author{   D.   Dacunha-Castelle$^1$  and   L.   Ferm\'{\i}n$^{1,2}$\\
{\scriptsize  $^1$ Universit\'e  Paris-Sud  \ and  \ $^2$  Universidad
Central            de            Venezuela}.\\            {\scriptsize
didier.dacunha-castelle@math.u-psud.fr,
fermin.lisandro@math.u-psud.fr} } 
\date{} \maketitle

\abstract{ We  prove that a large  set of long  memory (LM) processes
(including  classical LM  processes and  all processes  whose spectral
densities  have  a countable  number  of  singularities controlled  by
exponential  functions)  are  obtained  by an  aggregation  procedure
involving  short memory  (SM) processes  whose spectral  densities are
infinitely differentiable ($\mathcal{C}^{\infty}$).  We show that the
$\mathcal{C}^{\infty} $ class of spectral densities is the optimal
class to get a general result for disaggregation of LM processes in SM
processes,    in    the   sense    that    the    result   given    in
$\mathcal{C}^{\infty}$  class cannot be  improved taking  for instance
analytic  functions instead of  indefinitely derivable  functions.  \\
Keywords: Aggregation; disaggregation; long memory process; mixture.} 

\section{Introduction}

LM  processes, are  used in  many fields  such as  economics, finance,
hydrology or  communication networks.  Let $\gamma$  be the covariance
function of the process and
$$\|\gamma\|_{SM}=\sum_{k=0}^{\infty}|\gamma(k)|.$$
If $\|\gamma\|_{SM}<\infty$ then we say  that the process is SM and if
$\|\gamma\|_{SM}=\infty$,  we say that  the process  is LM.   The long
memory is  generally associated to  the singularities of  the spectral
density $f$.

Some of  these LM processes can be  seen as an aggregation,  that is a
normalized  sum  of  elementary  SM  processes,  see  \cite{Goncalvez,
Granger, Linden, Lippi&Zafferoni, Terence}.

All important  examples of LM  in finance, hydrology  or communication
networks etc.   are obtained by  a model of aggregation,  perhaps more
complicated that  the model  used here. Our  result says that  one has
almost  always the  possibility to  think  the LM  macroprocess as  an
aggregation   of  elementary  processes   with  $\mathcal{C}^{\infty}$
spectral density.

The  procedure  of aggregation  can  be  developed considering  doubly
random  elementary processes  $Z^i=\{Z^i_t(Y^i),\, t\in  T\}$ centered
second order and stationary, where $Y=\{Y^i, \, i \in \mathbb{N}\}$ is
a sequence  of i.i.d.  random  variables with distribution $\mu$  on $
\mathbb{R}^s$; $T=  \mathbb{Z}$ in the case of  discrete processes and
$T = \mathbb{R}$ in the continuous case, see \cite{Dacunha&Fermin1}.

For every fixed  trajectory of $Y$, we define  the sequence of partial
aggregations $X^N=\{X^N_t(Y),  \, t \in T\}$,  of elementary processes
$\{Z^i$\}, by
\begin{equation}\label{ec.1}
X^N_t(Y)=\frac{1}{B_N}\sum_{i=1}^N Z^i_t(Y^i),
\end{equation}
where  $B_N$ is a  normalization sequence  of positive  numbers. Under
some  general  conditions,  for  almost every  trajectory  $Y$,  $X^N$
converges in distribution  to the same Gaussian process  $X$, which is
called the aggregation of the elementary processes $\{Z^i\}$.

Equivalently, the  aggregation procedure  can be also  developed using
mixtures of  spectral densities as  main tool.  Let  $g(\lambda,y)$, $
\lambda \in I$, $ y \in \mathbb{R}^s$, be a spectral densities family,
where  $I=  (-\pi,  \pi]$  in  the  case  of  discrete  processes  and
$I=\mathbb{R}$  in  the continuous  case.   We  consider, for  $Y^i=y$
fixed, that  $g(\lambda,y)$ is the spectral density  of the elementary
process $Z^i(Y^i)$. Let $\mu$  be a probability on $\mathbb{R}^s$, the
$\mu$-mixture of the spectral densities $g(\lambda,y)$, is defined by
\begin{equation} \label{ec.2}
f(\lambda)=\int_{\mathbb{R}^s}g(\lambda,y)d\mu(y).
\end{equation}
$f(\lambda)$ is a well defined spectral density iff
\begin{equation}\label{ec.3}
\int_{I} f(\lambda)d\lambda < \infty.
\end{equation}

The  sequences of  processes  $\{X^N, N\in\mathbb{N}  \}$, defined  in
\eqref{ec.1},  is   now  the  partial  aggregation   sequence  of  the
elementary processes $\{Z^i\}$ associated to the mixture $f(\lambda)$.

Under  condition \eqref{ec.3}  in \cite{Dacunha&Fermin1}  it  is shown
that if the elementary  processes are independent then the aggregation
exists $\mu-a.s.$  and its spectral density  function is $f(\lambda)$.
Futhermore,  the  case  of  non independent  elementary  processes  is
studied  considering  an  interactive correlation  between  elementary
processes.

Disaggregation   is  the  inverse   procedure  of   aggregation.   Let
$\mathbb{G}$   be   a   class    of   spectral   densities   and   let
$\mathcal{M}(\mathbb{G})$  denote  the  set  of mixtures  of  spectral
densities  belonging to  $\mathbb{G}$.   We say  that  a process  with
spectral  density $f(\lambda)$  can be  disaggregated  into elementary
processes  with   spectral  densities  in  $\mathbb{G}$   iff  $f  \in
\mathcal{M}(\mathbb{G})$;    equivalently    by   disaggregation    in
$\mathbb{G}$, we  mean the  existence of a  representation of  a given
process  as  an  aggregation  of  elementary  processes  belonging  to
$\mathbb{G}$.  The  disaggregation problem  is then equivalent  to the
following question:  When do  we have $f  \in \mathcal{M}(\mathbb{G})$
for  $f$  and $\mathbb{G}$  given?   In  this  work, we  consider  the
procedure of aggregation only  for independent elementary processes. A
known  example is  the disaggregation  of the  $FARIMA(d)$  process on
$\mathbb{G}=\mathbb{AR}(1)$ the  class of autoregressive  processes of
order 1,  see \cite{ Lippi&Zafferoni}.  A more general  development of
the  disaggregation in  the class  of  $AR(p)$ processes  is given  in
\cite{Dacunha&Fermin2} where  we also show  results of disaggregation,
for continuous processes, in the class of $p$-order Orsteins-Uhlembeck
processes, $\mathbb{OU}(p)$.

Our purpose is  to study the disaggregation procedure  in SM processes
and on  subclasses $\mathbb{G}$ for which the  covariances decrease as
fast as  possible. We prove  that a large  set of LM  processes, whose
spectral densities have singularities of different kinds, are obtained
by  an  aggregation   procedure  involving  processes  whose  spectral
densities are infinitely differentiable.  There is no unicity for this
procedure. Classical LM processes are included in this set.

Conditions for the existence  of a disaggregation on $\mathbb{AR}(p)$,
$\mathbb{OU}(p)$  classes  or  on  $\mathcal{A}$,  set  of  analytical
spectral densities, are very restrictive.  They are related, as in the
$\mathbb{AR}(1)$ case,  to specific  algebraic properties of  $f$, for
instance to be a Mellin transform, and the elementary processes have a
spectral density roughly  of the form $\Phi (yf(\lambda  ))$ where the
regularity  of  $\Phi  $  is   sufficient  in  order  to  destroy  the
singularity of $f$.

We  show   that  the  class   $\mathcal{C}^{\infty}$  of  indefinitely
differentiable spectral densities is the  right class to get a general
result for disaggregation of LM  processes in SM processes.  For $f\in
$ $\mathcal{C}^{\infty}$, the respective covariance function decreases
faster than $n^{-j}$ for any $j$.  This decay can measure the level of
SM  comparatively  to  that  of $\mathcal{A}$  which  is  exponential.
Disaggregation  on  $\mathcal{C}^{\infty}$  is  linked  to  very  weak
analytical properties much easier to satisfy by $f$ than the algebraic
properties  required  by analytic  disaggregation  and the  elementary
densities    are   here   roughly    of   the    form   $g(\lambda,y)=
f(\lambda)K(\lambda,y)$,  where $K$  is  a kernel  such  that at  each
singularity $\lambda_0$ of  $f$ we have $f(\lambda_0)K(\lambda_0,y)=0$
as  well  as for  all  its  derivatives.  The $\mathcal{C}^{\infty  }$
behavior  of  $K$  drives  the  $\mathcal{C}^{\infty  }$  behavior  of
$g(\lambda,y)$.   Then  we  show  examples  of  $\mathcal{C}^{\infty}$
functions that can not be disaggregated on the class $\mathcal{A}$. On
the  other hand,  for very  rough densities,  we give  examples  of LM
processes which  cannot be the aggregation of  short memory processes.
We shall call  such a situation as ''hard''  long memory.  This result
is  a  good  illustration  of   ''how  much  singularity  have  to  be
concentrated'' in order to generate LM by aggregation.

\section{Condition  for the Existence  of Disaggregation
in SM}

Let us begin  with some remarks about mixtures  and LM property giving
rough  necessary conditions for  existence of  a disaggregation  on SM
densities. We use  only the fact that every  SM density is continuous,
(in fact sufficient conditions for  SM property are for instance to be
$\alpha$-lipchitzian  for every  $\alpha, \,  0<\alpha <1$,  or  to be
$\alpha$-lipchitzian for  some $\alpha  > \frac{1}{2}$ and  of bounded
variation).

\begin{lemma}\label{lemma1}
A spectral density  is a mixture of continuous  spectral densities iff
it is lower semicontinuous (l.s.c.).
\end{lemma}
\begin{proof}
If $f$  is such a  mixture, applying Fatou's  lemma we see that  it is
l.s.c.   Conversely, if  $f$ is  l.s.c.   and positive,  it exists,  a
sequence    $g_{n}$   of   continuous    functions   such    that   $\
f=\sup_{n}g_{n}$.  Taking the supremum until $n$ we can choose $g_{n}$
increasing.  So  $f=\sum_{n\geq 1} 2^{-n}(2^{n}(g_{n+1}-g_{n}))$ which
is of the form $\int g d\mu $.
\end{proof}

From the lemma, we see that  every non l.s.c.  spectral density has LM
property and cannot be disaggregated on SM class. For instance, if $h$
is  the function  equal to  one on  a perfect  set,  \cite{Kahane}, of
strictly positive  Lebesgue measure, then $h$  is upper semicontinuous
(u.s.c.)   but non  l.s.c.   and it  is  the density  of a  absolutely
continuous  probability with  respect to  the Lebesgue  measure.  This
provides  an example  of  a situation  that  we can  call ''hard''  LM
process which  cannot be disaggregated  by SM processes.  In  fact, we
are very  closed to the case  of a non  absolutely continuous spectral
measure.

 \section{Disaggregation     of    LM    Processes    on
$\mathcal{C}^{\infty}$ Class}

Long  memory of  a process  with spectral  density $f$  is  in general
associated to singularities of $f$ or  of some of its derivatives at a
frequency $\lambda_0$.  Singularities are  often classified as a first
order when one side limit exists and second order singularity when the
function tends to  $\infty $ as $\lambda $ tends  to $\lambda_0$ or it
does not exist any limit at $\lambda_0$ (the function being bounded or
not, with bounded variation or not).  We try to take into account most
of these  situations. Our  main purpose  is to obtain  for a  class as
broad as possible, including all classical examples but not limited to
more  or  less  explicit  densities, a  disaggregation  on  elementary
processes with the  best possible decay of correlations.   We are lead
to  work  mainly   with  $\mathbb{G}  =  \mathcal{C}^{\infty}$  (resp.
$\mathcal{C}^{H}$) the  class of  all spectral densities  belonging to
$\mathcal{C}^{\infty  }$  (resp.   $\mathcal{C}^{H}$);  an  equivalent
property   is    that   the    correlation   function   of    $f   \in
\mathcal{C}^{\infty}$ tend  to $0$  faster than ${1}/{n^{j}}$  for any
$j\in   \mathbb{N}$  (resp.   faster   than  ${c}/{n^{H}}$   for  some
$c>0$). For  $H$, $1\leq H \leq  \infty$, given we get  for a function
$f\in  \mathcal{C}^{H}$  except  for  a  finite or  countable  set  of
frequencies, a  disaggregation on  a $\mathcal{C}^H$ class.   It seems
that disaggregation  on the class $\mathcal{C}^{\infty}$  is easier to
reach than a disaggregation on the SM class.

In order  to extend  the qualitative situation  analyzed for  the case
$\mathbb{G}=\mathbb{AR}(1)$,  in \cite{Dacunha&Fermin2},  we introduce
the following definition.

\begin{definition} \label{def.1}
  Let $\Lambda  \subset I$ be a  finite set of  frequencies, we define
  $F_{\Lambda}^{H}$   as   the    set   of   spectral   densities   in
  $\mathcal{C}^H$, $1 \leq H \leq \infty$, which are 0 on $\Lambda$ as
  all   their  $H$   derivatives.   If   $\Lambda=\{0\}$   we  denote
  $F_{\Lambda}^{H}$ by $F_0^{H}$.
\end{definition}

Let now $K$  be a positive kernel, for  instance on $\mathbb{R}^+$, so
$\int_0^{\infty}K(y)dy=1$  and let  $\phi$ be  a positive  function on
$I$. We use the trivial relation
\begin{equation}\label{ec.4}
f(\lambda )=\int_0^{\infty}g(\lambda ,y)dy
\end{equation}
where $g(\lambda,y) = f(\lambda)\phi(\lambda)K(\phi(\lambda)y)$.

The  pair $[K;\phi]$  will be  called a  killer kernel  (it  kills the
singularities).   The  generic  situation  is  given  by  the  formula
$f(0)\phi(0)K(\phi(0)y)=\lim_{\lambda          \rightarrow          0}
f(\lambda)\phi(\lambda)K(\phi(\lambda)y)=0$    for   every    $y   \in
\mathbb{R}^+$,  even if  $f(0)=\infty$. So  f is  the $dy$  mixture of
$g(\lambda,y) \in F_0^H$.  We can extend the definition of mixtures of
spectral  densities, given  in \eqref{ec.4},  taking for  $\mu$  a non
bounded measure,  as the  Lebesgue measure $dy$,  in this case  we can
consider  that we  take  a strictly  positive  density of  probability
$\sigma(y)$    and    then     we    consider    the    mixtures    of
$\sigma^{-1}(y)g(\lambda,y)$ by the probability $\sigma(y)dy.$

Let us give  some examples of disaggregation of a  function $f$ with a
single singularity at $\lambda_0$.  We consider the standard situation
where $f$ and its derivatives are explicitly controlled by exponential
finctions. The key  to solve our problem is  the representation of $f$
given in \eqref{ec.4}.  
\\
\textbf{Example    1.}     Let
$f(\lambda)=1_{(-\lambda_0,\lambda_0)}(\lambda)$  for  $\lambda_0  \in
I$,  $\phi(\lambda)=  1/|\lambda^2-\lambda_0^2|^p$  with $0<p<1$,  and
$K(y)=e^{-y}$.  In  this case is  easy to check  that $g(\lambda,y)\in
F^{\infty}_{\{-\lambda_0, \lambda_0\}}$, since  all derivatives of $g$
are 0 for $|\lambda|=\lambda_0$. 
\\
\textbf{Example  2.}   Let
$f(\lambda)=|1-\cos(  \lambda -  \lambda_0)|^{-d}$,  $0<d<1$, $\lambda
\in I=(-\pi,\pi]$.  We keep $\phi(\lambda)= 1/|\lambda - \lambda_0|^p$
with   $0<p<1$,  and   $K(y)=e^{-y}$.  All   derivatives  of   $f$  at
$\lambda=\lambda_0$   are   controlled   by   a  negative   power   of
$|\lambda-\lambda_0|$ and  so $g(\lambda,y) \in  F_{\lambda_0}^H$. The
same  properties can  be  easily checked  for  a strongly  oscillating
function $f$ as $\cos(\pi (\lambda-\lambda_0))/|\lambda-\lambda_0|^q$,
for $0<q<1$.  So for these  kinds of controlled singularities, we show
that  $f\in  \mathcal{M}(F^H_{\Lambda}) \subset  \mathcal{M}(C^{\infty
}).$

\begin{definition} \label{def.2}
Let $f$  be a spectral  density. We say  that $f \in S^{H}$,  $1\leq H
\leq  \infty$,  iff $f$  has  a  continuous  $H$ derivative  at  every
frequency except for a finite  set of frequencies $\Lambda = \{\lambda
_{j}, j\in J\}$  and if it  exists $q$, $0<q<1$,  and $a$, $0<a<1$,
such that for all $j\in J$ and for all $l\leq H$
$$\lim_{\lambda \rightarrow  \lambda_j}\exp\left( -\frac{a}{|\lambda -
\lambda_j|^q}\right)|f^{(l)}(\lambda)| = 0.$$
If $\Lambda$  is a  countable infinite set  instead of finite  and has
only a finite  number of accumulation points, then we  say that $f \in
T^{H}$.
\end{definition}

We state now a theorem for the standard situation.

\begin{theorem} \label{theo.1}
Let    $f\in    S^{H}$,    $1\leq    H\leq   \infty$,    then    $f\in
\mathcal{M}(F^H_{\Lambda}) \subset \mathcal{M}(\mathcal{C}^{H})$.
\end{theorem}

\begin{proof}
Let                                    $\phi(\lambda)=\left(\prod_{j\in
J}\left|\lambda-\lambda_{j}\right|^p\right)^{-1}$,  $K(y)=e^{-y}$  and
$g(\lambda,   y)=   f(\lambda)\phi(\lambda)K(\phi(\lambda)y)$.    Then
$f(\lambda )=\int_0^{\infty}  g(\lambda ,y)dy$. We  choose  $p$ such that
$0<q<p<1$. If  $\Psi(\lambda, y)= \phi(\lambda)K(\phi(\lambda)y)$ then
for  all the $l$-derivatives  of $g$,  $l\leq H$,  we have  that it exists
constants $b_l$, $C_l$ and $m_l$ such that
\begin{equation*}
\left|g^{(l)}(\lambda,y)\right|
=\left|\sum_{k=0}^lC_{k,l}f^{(k)}\psi^{(l-k)}(\lambda,y)\right|    \leq
C_l|\Psi(\lambda,y)|^{m_l}e^{\left(-\sum_{j\in                J}\frac{a
b_l}{|\lambda-\lambda_j|^q} \right)}.
\end{equation*}
So    $g(\lambda,     y)    \in    F^H_{\Lambda}$     and    $f    \in
\mathcal{M}(F^H_{\Lambda}) \subset \mathcal{M}(C^{H})$.
\end{proof}
\qquad\\
\textbf{Example  3.}    If   $f(\lambda  )=|\lambda|^{-d}$,
$\lambda \in  I=\mathbb{R}$, $0<d<1$, $f$  is the spectral  density of
continuous   fractional  Brownian   motion.    Then  $|f^{(l)}(\lambda
)|=O(1/|\lambda|^{d+l})$  and the  conditions of  Theorem \ref{theo.1}
are satisfied.

Next  theorem is  an extension  of  Theorem \ref{theo.1}  but using  a
localization required if $\Lambda $ is infinite.

\begin{theorem}\label{theo.2}
Theorem \ref{theo.1} remains valid if $f\in T^{H}$.
\end{theorem}

\begin{proof}
Suppose,  in  order to  simplify  notations,  that  $\Lambda $  is  an
infinite countable  set with only one accumulation  point. The general
case can be easily obtained by re-indexing of $\Lambda$ points using the
partition of  $I$ defined by  the points of accumulation  of $\Lambda$
and   then   applying   the   same   proof.    So   we   can   suppose
$\Lambda=\{\lambda_j,  j \geq 1\}$,  with $\lambda_j  < \lambda_{j+1}$
for every $j \in \mathbb{N}$.

We  build in the  same spirit  than previously  a family  of functions
$g^{\Lambda}(\lambda,  y)$,   multiplying  $f$  by   a  killer  kernel
$[K;\phi^{\Lambda}](\lambda,  y)$  that   annihilates  the  points  of
discontinuity of $f$.

Let us note $a=\inf I$,  $b= \sup I$  and $\lambda_{\infty}=\lim_{j
 \rightarrow \infty} \lambda_j =  \sup_j \lambda_j$. Let be $p$ such that
 $0 < q < p < 1$, and we consider
\begin{eqnarray*}
\, [K, \phi^{0}](\lambda,y) & = & \frac{1}{|\lambda-\lambda_1|^p} \exp
\left(-\frac{y}{|\lambda-\lambda_1|^p}\right)1_{(a,\lambda_1)}(\lambda).
\\ \, [ K, \phi^{j}](\lambda,y) & = & \frac{1}{|\lambda - \lambda_j|^p
|\lambda   -  \lambda_{j+1}|^p}   \exp  \left(-   \frac{y}{|\lambda  -
\lambda_j|^p   |\lambda   -  \lambda_{j+1}|^p}\right)   1_{(\lambda_j,
\lambda_{j+1})}(\lambda).  \\ \, [  K,\phi^{\infty}] (\lambda,y) & = &
\frac{1}{|\lambda-\lambda_{\infty}|^p}\exp  \left(-\frac{y}{|\lambda -
\lambda_{\infty}|^p}\right) 1_{(\lambda_{\infty},b)}(\lambda).
\end{eqnarray*}
Then we define
$$[K;\phi^{\Lambda}](\lambda,y)=\sum_{j=0}^{\infty}[K;\phi^j](\lambda,y),$$
and
$$g^{\Lambda}(\lambda,y)=f(\lambda)[K;\phi^{\Lambda}](\lambda,y).$$
We have that,
$$\int_0^{\infty}[K;\phi^{\Lambda}](\lambda,y)dy=\int_0^{\infty}e^{-z}dz
\left(1_{(a,\lambda_1)}+\sum_{j                 \geq                1}
1_{(\lambda_j,\lambda_{j+1})}+1_{(\lambda_{\infty},b)}\right) = 1. $$
So  $V(y)=\int_{I}g^{\Lambda }(\lambda ,y)d\lambda  $ and  by applying
Fubini's  theorem $\int_{0}^{\infty}V(y)dy=\int_{I}f(\lambda)d\lambda=
\gamma_0<\infty$.

We can  prove, by  using the same  proof as for  Theorem \ref{theo.1},
that the  $H$ derivatives  with respect to  $\lambda $  of $g^{\Lambda
}(\lambda,y)$ go to $0$  when $\lambda \rightarrow \lambda_{j}$, since
$q<p$.  So $g^{\Lambda }(\lambda,y) \in F^{H}_{\Lambda}.$
\end{proof}

\begin{remark}
Killer  kernels $[K;\phi]$  selected to  build only  the  mixtures are
never the  best ones,  for coveriances decay.   For instance,  we can
take  $\exp(-\exp(y))$  instead   of  $\exp(-y)$  getting  covariances
decreasing to $0$ slightly faster and so on.
\end{remark}

\begin{remark}
We can define mixtures of  Wold regular densities, in the Wold Theorem
sense, \cite{Azencott},   which  verify  the
following condition
\begin{equation}\label{ec.5}
\int_{I} \log g(\lambda,y) d\lambda > -\infty, \quad \mu -a.s.
\end {equation}
In  this case,  we  say  that the  processes  with spectral  densities
$g(\lambda,y)$ are  regular.  If $f(\lambda)$  is the mixture  of the
densities $g(\lambda,y)$,  condition \eqref{ec.5} does  not imply that
$\int_{I} \log f(\lambda)  d\lambda > -\infty$. But if  $f$ is regular
then $g(\lambda,y)$ is regular $\mu -a.s.$ by Jensen's inequality. The
main  point on  this topic  is that  we can  choose the  killer kernel
$[K;\phi]$  such that if  $f$ is  Wold regular  then all  the elementary
processes used in the aggregation are also Wold regular.
\end{remark}

\section{Analytic spectral densities}

Let us prove, in a certain  sense, that the previous results cannot be
improved taking  analytic functions instead  of indefinitely derivable
functions.

Disaggregation    is    a    hierarchical    procedure:    if    $f\in
\mathcal{M}(\mathbb{G})$         and        $\mathbb{G}        \subset
\mathcal{M}(\mathbb{H})$ then  $f\in \mathcal{M}(\mathbb{H})$, in fact
if $g(\lambda,y)=\int h(\lambda,z)d\nu(y,z)$ then
$$f(\lambda   )=\int   g(\lambda   ,y)d\mu(y)=\int   \int   [h(\lambda
    ,y,z)d\nu(y,z)]d\mu (y).$$
In general we have $ \mathbb{G} \subset \mathcal{M}(\mathbb{G})\subset
\overline{\mathcal{M}(\mathbb{G})}$ with strict inclusion, the closure
being  taken  in  $L^{1}(d\lambda)$.    The  obvious  exception  is  $
\mathbb{G} =\mathbb{MA}(q)$, the  set of $q$-moving average densities,
for which $  \mathbb{MA}(q) = \mathcal{M}(\mathbb{MA}(q)) = \overline{
\mathcal{ M} (\mathbb{MA}(q))}.$

We use  this hierarchical procedure in  order to show  that our result
can not  be improved in the  following sense: we  cannot take analytic
functions instead of $\mathcal{C}^{\infty}$.  So we have to check that
the   functions  we   have   used  in   $\mathcal{C}^{\infty  }$,   as
$f(\lambda)=\frac{1}{|\lambda|^p}\exp(-\frac{y}{|\lambda   |^p})$,  do
not belong  to $\mathcal{M}(\mathcal{A})$, in  order to show  that our
result can not be improved.

Let      $f      \in      \mathcal{C}_{0,+}^{\infty      }$      where
$\mathcal{C}_{0,+}^{\infty   }=\left\{  f:   \,  f\geq   0,   \,  f\in
\mathcal{C}^{\infty},  \text{ it  exists} \,  \lambda_0  \, \text{such
that}\,   f^{(j)}(\lambda_0)=0    \,   \text{for   every}    \,   j\in
\mathbb{N}\right\}$.                 Suppose                $f(\lambda
)=\int_{\mathbb{R}^{s}}g(\lambda, y)d\mu  (y)$ with $g\in \mathcal{A}$
$\mu-a.s.$ The Fatou's Lemma implies
$$f^{(j)}(\lambda_0)\geq        \int_{\mathbb{R}^{s}}g^{(j)}(\lambda_0,
y)d\mu (y),$$
and  $f^{(j)}(\lambda  _0)= 0$  implies  that, if  $g^{(j)}(\lambda_0,
y)\geq 0$ $\mu-a.s.$, then  $ g^{(j)}(\lambda_0, y)=0$ $\mu-a.s.$ From
$g(\lambda  ,y)\geq   0$  $  \mu-a.s$  we   get  $\  g(\lambda_0,y)=0$
$\mu-a.s.$  and $g^{(1)}(\lambda_0,y)  \geq  0$ $\mu-a.s.$,  and so  $
g^{(1)}(\lambda_0, y)= 0$ $\mu -a.s$ and $g^{(2)}(\lambda _{0}, y)\geq
0$ $\mu -a.s$. By induction we have that $g^{(j)}(\lambda_0,y)=0$ $\mu
-a.s$ and  $\mu \{y,  g(\lambda ,y)\in \mathcal{A}\}  = 0 $.   We have
proved that $f(\lambda) \notin \mathcal{M}(\mathcal{A})$.

Disaggregation in  the class of  analytic functions is a  very limited
possibility. A slight modification of  the proof just above shows that
spectral  densities which  are constant  or linear  or  polynomial (of
given degree) by pieces cannot be in $\mathcal{M}(\mathcal{A})$ except
if they are themselves elements of $\mathcal{A}$, that is, if they are
polynomials.

\qquad\\ \textbf{Acknowledgements:} The research of L. Ferm\'{\i}n was
supported  in part by  a grant  from the  FONACIT and  Proyecto Agenda
Petr\'oleo (Venezuela).  \\

\end{document}